\documentclass[11pt]{article}

\usepackage{latexsym}
\usepackage{amsfonts}
\usepackage{amssymb}
\usepackage{amsmath}
\usepackage{mathrsfs}
\usepackage[all]{xy}
\usepackage{graphicx}
\usepackage{stmaryrd}

\newtheorem{Theory}{Theory}[section] %Counter for all types of Theorems
\newtheorem{theorem}[Theory]{Theorem}
\newtheorem{lemma}[Theory]{Lemma}

\newtheorem{corollary}[Theory]{Corollary}

\newtheorem{remark}[Theory]{Remark} %trivial but worth noticing
\newtheorem{claim}[Theory]{Claim}

\newcommand{\kab}{K_{(a,b)}}
\newcommand{\fk}[2]{K_{(#1,#2)}}
\newcommand{\tk}[2]{\tilde{K}_{(#1,#2)}}
\newcommand{\be}{\begin{enumerate}}
\newcommand{\ee}{\end{enumerate}}
\newcommand{\bt}{\begin{theorem}}
\newcommand{\et}{\end{theorem}}
\newcommand{\bl}{\begin{lemma}}
\newcommand{\el}{\end{lemma}}
\newcommand{\bc}{\begin{corollary}}
\newcommand{\ec}{\end{corollary}}
\newcommand{\br}{\begin{remark}}
\newcommand{\er}{\end{remark}}
\newcommand{\bcl}{\begin{claim}}
\newcommand{\ecl}{\end{claim}}
\newcommand{\inner}{\operatorname{Inner}}
\renewcommand{\outer}{\operatorname{Outer}}
\newcommand{\rev}{\operatorname{Rev}}
\newcommand{\res}{\operatorname{Res}}
\newcommand{\supp}{\operatorname{Supp}}

\begin{document}

\begin{center}
{\LARGE{{\bf{Finite Index Subgroups of R. Thompson's Group F}}}}

\vspace{10pt}
Collin Bleak and Bronlyn Wassink
\end{center}

\vspace{10pt} {\flushleft ABSTACT: The authors classify the finite
index subgroups of R. Thompson's group $F$.  All such groups that are
not isomorphic to $F$ are non-split extensions of finite cyclic groups
by $F$.  The classification describes precisely which finite index
subgroups of $F$ are isomorphic to $F$, and also separates the
isomorphism classes of the finite index subgroups of $F$ which are not
isomorphic to $F$ from each other; characterizing the structure of the
extensions using the structure of the finite index subgroups of
$Z\times Z$.}

\vspace{15pt}
\section{Introduction}
In this paper, we classify the finite index subgroups of R. Thompson's
group $F$.  By this classification, we are able to answer Victor
Guba's Question 4.5 in the problems report \cite{Thompson40}.

The group $F$ was introduced by R. Thompson in the late 1960's as part
of a family of groups $F\leq T\leq V$.  It has been the object of much
study, and it's theory has impacted various fields of mathematics,
including not only the theory of infinite groups, but also
low-dimensional topology, simple homotopy theory, measure theory, and
even category theory.  An introductory reference to the
theory of $F$, $T$, and $V$ is the survey paper \cite{CFP}.

The main characterization of $F$ that we will use is that it is the
group of all piecewise-linear, orientation-preserving homeomorphisms
of the unit interval which admit finitely many breaks in slope, where
these breaks are restricted to occur over the diadic rationals
$Z[1/2]$, and where all slopes of affine segments of the graphs of
these elements are integral powers of two.  We will also use some of
the standard presentations for $F$ in our analysis, which
presentations will be given in the next section.

We will use the notation $FIF$ to represent the set of all finite index subgroups of $F$.  

In order to state our results in full, we will need to build a
specific homomorphism.

Given $f\in F$, we will denote the derivative of
$f$ at $x$ by $f'(x)$, if it exists.  We will also define $f'(0)$
to be the derivative from the right at $0$, and $f'(1)$ to be the
derivative from the left at $1$.  Note that these last two
derivatives always exist as elements of $F$ are affine near $0$ and
$1$ in $(0,1)$.

We now define a well known homomorphism $\phi:F\to Z^2$ by the rule
\[
\phi(f) = (\log_2(f'(0)), \log_2(f'(1)))
\] 
for all $f\in F$.

By a standard fact in the literature of $F$ (Theorem 4.1 from
\cite{CFP}), the group commutator subgroup $F'$ of $F$ consists of
precisely the elements in $F$ with leading and trailing slopes one, that is,
$F'$ is the kernel of the map $\phi$.

Given two positive integers $a$ and $b$, we define $\tilde{K}_{(a,b)}
= \langle(a,0),(0,b)\rangle\leq Z^2$.  We now define 
\[
\kab =
\phi^{-1}(\tilde{K}_{(a,b)}).
\]
In particular, $\kab$ can be thought of as the group of all elements
in $F$ with graphs having slopes near zero as integral powers of $2^a$
while at the same time having slopes near one as integral powers of
$2^b$.  We will call any such $\kab\leq F$ a \emph{rectangular
subgroup of $F$}, or simply a \emph{rectangular group}.  We will also refer to the groups $\tk{a}{b}\leq Z^2$ as rectangular groups, where the context will make clear which sort of rectangular groups we are refering to.

We are now ready to give an explicit list of our results.  Our first
theorem is a corollary of our last theorem, but as we will prove it
earlier in the paper in a direct fashion, we will list it here as a
stand-alone result.

\bt\label{KabIsoF} 

Let $H\in FIF$.  $H$ is isomorphic to $F$
if and only if $H = \kab$ for some positive integers $a$ and $b$.

\et

Given positive integers $a$ and $b$, it is immediate that $F/\kab\cong
Z_a\times Z_b$, in particular we have the following theorem.

\bt
\label{FIsExtension} 

Given any positive integers $a$ and $b$, $F$ can be regarded as a
non-split extension of $Z_a\times Z_b$ by $F$.  In particular,
there are maps $\iota$ and $\tau$ so that the following sequence is
exact.

\[
\xymatrix{
1\ar[r]&F\ar[r]^{\iota}&F\ar[r]^{\!\!\!\!\!\!\!\!\!\!\!\!\tau} & 
Z_a\times Z_b\ar[r] &1.}
\]

\et

We will prove two further theorems.  Before stating them, we mention
some key lemmas, and build some language that will help with the
statements of the theorems.

\bl \label{FPrimeInFIF}
If $H\in FIF$ then $F'\leq H$.
\el

Once the previous lemma is established, it is not hard to come to the
following lemma.

\bl\label{FIFNormal}
If $H\in FIF$ then $H\trianglelefteqslant F$.
\el

Now, the above lemmas assure us that we can analyze all of the finite
index subgroups of $F$ by considering the finite index subgroups of $Z^2$.

We have one further lemma, which will assist us in our statements below.

\bl \label{innerOuter}

Suppose $H\in FIF$. There exist rectangular groups $\inner(H)$ and $\outer(H)$ so that $\inner(H)$ is a unique maximal
rectangular subgroup in $H$ and $\outer(H)$ a unique minimal
rectangular group containing $H$.

\el
In particular, if $H\in FIF$, then we have the following list of containments (where the first two are equalities in the case that $H$ is a rectangular group).
\[
\inner(H)\trianglelefteqslant H\trianglelefteqslant\outer(H)\trianglelefteqslant F
\]

We are now in a position to state our next theorem.

\bt\label{FIFCorreFIZ2} 
\be
\item The map $\phi$ induces a one-one correspondence between the
finite index subgroups of $F$ and the finite index subgroups of $Z^2$.
\item Suppose $H$ is a finite index subgroup of $F$, with image
$\tilde{H}=\phi(H) \leq Z^2$, and that $a$ and $b$ are positive
integers so that $\inner{H}=\kab\leq H$.  If
$Q=\tilde{H}\slash\tilde{K}_{(a,b)}$, then $Q$ is finite cyclic, and
there are maps $\iota$, $\rho$, $\tilde{\iota}$ and $\tilde{\rho}$ so
that the diagram below commutes with the two rows being exact:
\[
\xymatrix { 
&  F\ar[d]_{\cong}  &  &  &  \\ 
1\ar[r]  &  {\kab}\ar[r]^{\iota} \ar[d]_{\phi|_{\kab}}  & 
H\ar[r]^{\rho}\ar[d]_{\phi|_H} & Q\ar@{=}[d]\ar[r] & 1\\ 
1\ar[r]  &  {\tk{a}{b}}\ar[r]^{\tilde{\iota}}  & 
{\tilde{H}}\ar[r]^{\tilde{\rho}}  &  Q\ar[r]  &  1.
}
\]
\ee
\et

The essence of the above theorem is that in $Z^2$, each finite index
subgroup $\tilde{H}$ is a finite cyclic extension (by $Q$ above) of the maximal
rectangular subgroup of $\tilde{H}$.  The extension pulls back, so that the
finite index subgroup $H$ of $F$ can be seen as a finite cyclic
extension of the maximal rectangular group $\kab$ in H by the same
group $Q$ .  Whenever $Q$ is non-trivial, the resulting extension is
non-split and results in a group that is not isomorphic with $F$.  

We will give several examples at the end of the paper where $Q$ above
is non-trivial, that is, examples of finite index subgroups of $F$
which are not isomorphic to $F$.

We
define Res:$\,FIF\to N$, where we use the rule $H\mapsto n$, where $n$
is the cardinality of $H/\inner(H)$.  We will call the value $n$ in
the last sentence the \emph{residue} of $H$.

It turns out the relationship between a finite index subgroup of $F$
and its maximal rectangular subgroup is very special.  We show the
following lemma.

\bl\label{characteristicKab} 

Suppose $H$, $H'\in FIF$, $K=\inner(H)$, $K' = \inner(H')$, and
$\xi:H\to H'$ is an isomorphism. Then

\be
\item $\xi(K)=K'$
\item $K$ is characteristic in $H$ and $K'$ is characteristic in $H'$, and
\item $\res(H) = \res(H')$.
\ee
\el

Note that in the above, the second two points follow easily from the first.

For convenience, given $a$, $b$ positive integers, let us fix a
particular isomorphism $\tau_{(a,b)}\!:\kab\to F$, so that if
$f\in\kab$ with $\phi(f) = (as,bt)$ then $\phi(\tau_{(a,b)}(f)) = (s,t)$.  (We
note that these are the precise sorts of isomorphisms which we build
in the proof of Theorem \ref{KabIsoF}.)  We also need to name the
isomorphism Rev:$\,F\to F$ which is obtained if we conjugate the
elements of $F$ by the orientation-reversing map rev$:[0,1]\to[0,1]$
defined by the equation rev$(x) = 1-x$.  We are now ready to state our
final theorem.

\bt \label{isoClassification} 

Suppose $H$, $H'$ are finite index subgroups of $F$.  Let $a$, $b$,
$c$, $d$ be positive integers so that $\kab= \outer(H)$,
$\fk{c}{d}=\outer(H')$. $H$ is isomorphic with $H'$ if and only if
$\tau_{(a,b)}(H) = \tau_{(c,d)}(H')$ or $\tau_{(a,b)}(H) =
\rev(\tau_{(c,d)}(H'))$.  
\et

These investigations were started when Jim Belk asked the first author
if he knew whether or not [the group we call $K_{(2,2)}$] is
isomorphic to $F$.  The approach taken in this paper was motivated by
the proof of Brin's ubiquity result (see \cite{BrinU}), where Brin
shows that a subgroup of the full group of piecewise-linear,
orientation-preserving homeomorphisms of $[0,1]$ contains a copy of
R. Thompson's group $F$ if certain weak geometric conditions are
satisfied.

We are unaware of any published results relating to our own work here.
However, Burillo, Cleary and R\"{o}ver, in the course of their investigations
into the abstract commensurator of $F$, and using techniques different from
our own, have also understood the one-one correspondence between the finite 
index subgroups of $Z^2$ and the finite index subgroups of $F$.  Also, they 
have the result that the rectangular subgroups of $F$ are isomorphic to $F$. 
See \cite{BCRCommensurator}.

The authors would like to thank Matt Brin for interesting
discussions of these results, and also for some observations and
questions which helped us to refine the results.  Also, the first
author would like to thank Jim Belk for asking the initial question
that lead to this work, and to thank Mark Brittenham, Ken Brown, Ross
Geoghegan, Susan Hermiller, and John Meakin for interesting
conversations about these results.

\section{Definitions and Notation}

Richard Thompson's Group $F$ can also described by the following
presentations.

\[ F \cong \langle x_0, x_1, x_2, ... \ | \ x_j^{x_i} = x_{j+1} \ \mbox{for} \ i<j \rangle \]
\[ F \cong \langle x_0, x_1 \ | \ [x_0x_1^{-1}, x_1^{x_0}]=[x_0x_1^{-1}, x_1^{x_0^2}]=1   \rangle\] 
where $a^{b} = b^{-1}ab$ and $[a,b] = aba^{-1}b^{-1}$.

In these presentations, the generators $x_0$ and $x_1$ can be realized
as piecewise-linear homeomorphisms of the unit interval with breaks in
slope occurring over the diadic rationals, and with all slopes being
integral powers of two (that is, as elements of $F$ using the
definition of $F$ as a group of homeomorphisms of the unit interval).
We establish the mechanism of specifying any such function by listing
the points in its graph where slope changes.  We will call such points
\emph{breaks}, so that we will specify an element of $F$ by listing
its set of breaks.

Let $f_0$ be the element with breaks
$\left\{(1/4,1/2),\,(1/2,3/4)\right\}$ and let $f_1$ be the element
with breaks $\left\{(1/2,1/2),\,(5/8,3/4),\,(3/4,7/8)\right\}$.  The
functions $f_0$ and $f_1$ play the roles of $x_0$ and $x_1$ in the
presentations above.  Here are the graphs of these functions.
\begin{center}
\includegraphics[scale=0.25]{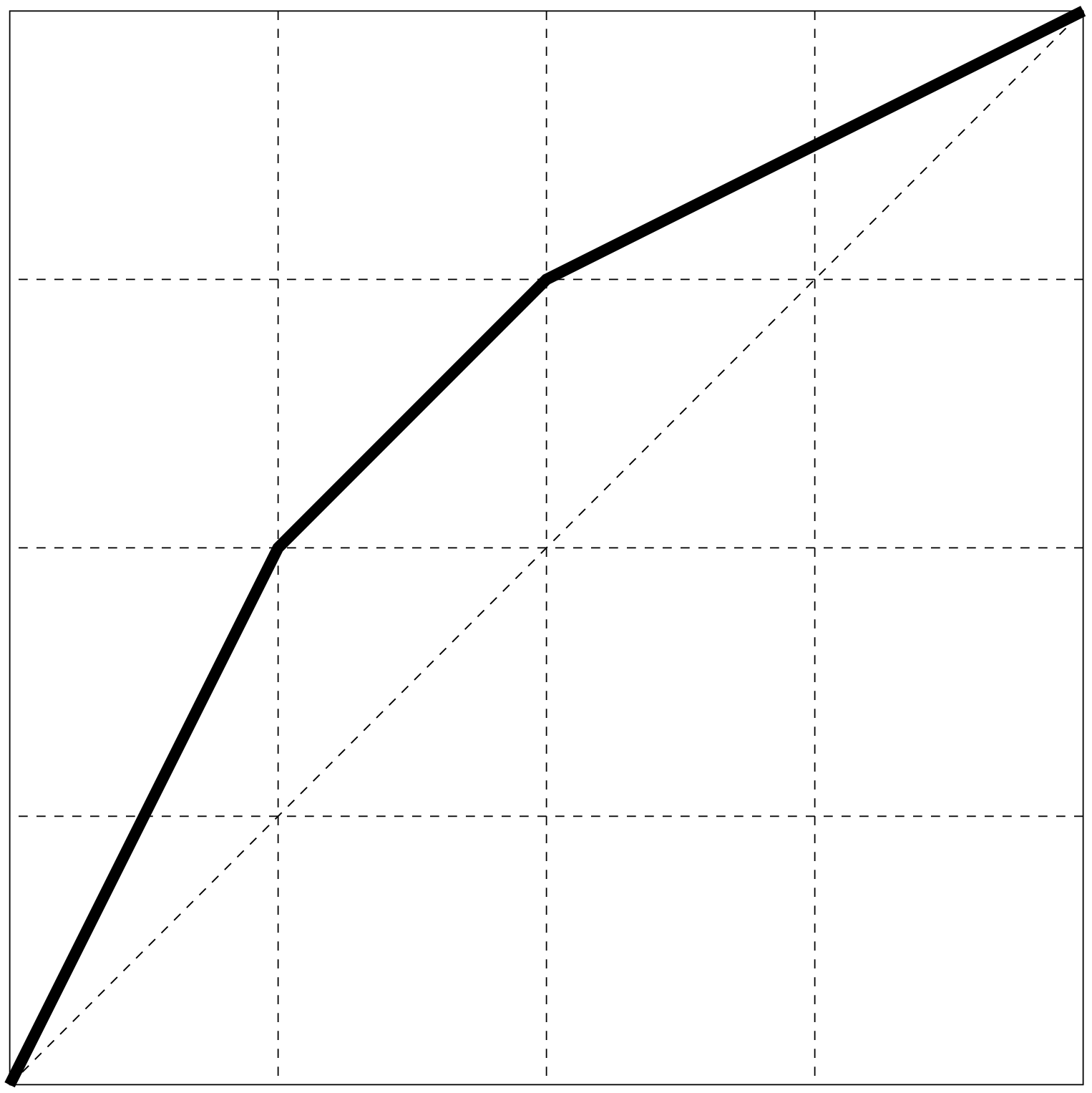} \ \ \ \ \ \ \ \ \ \ \ \ \ \ \ \ 
\includegraphics[scale=0.25]{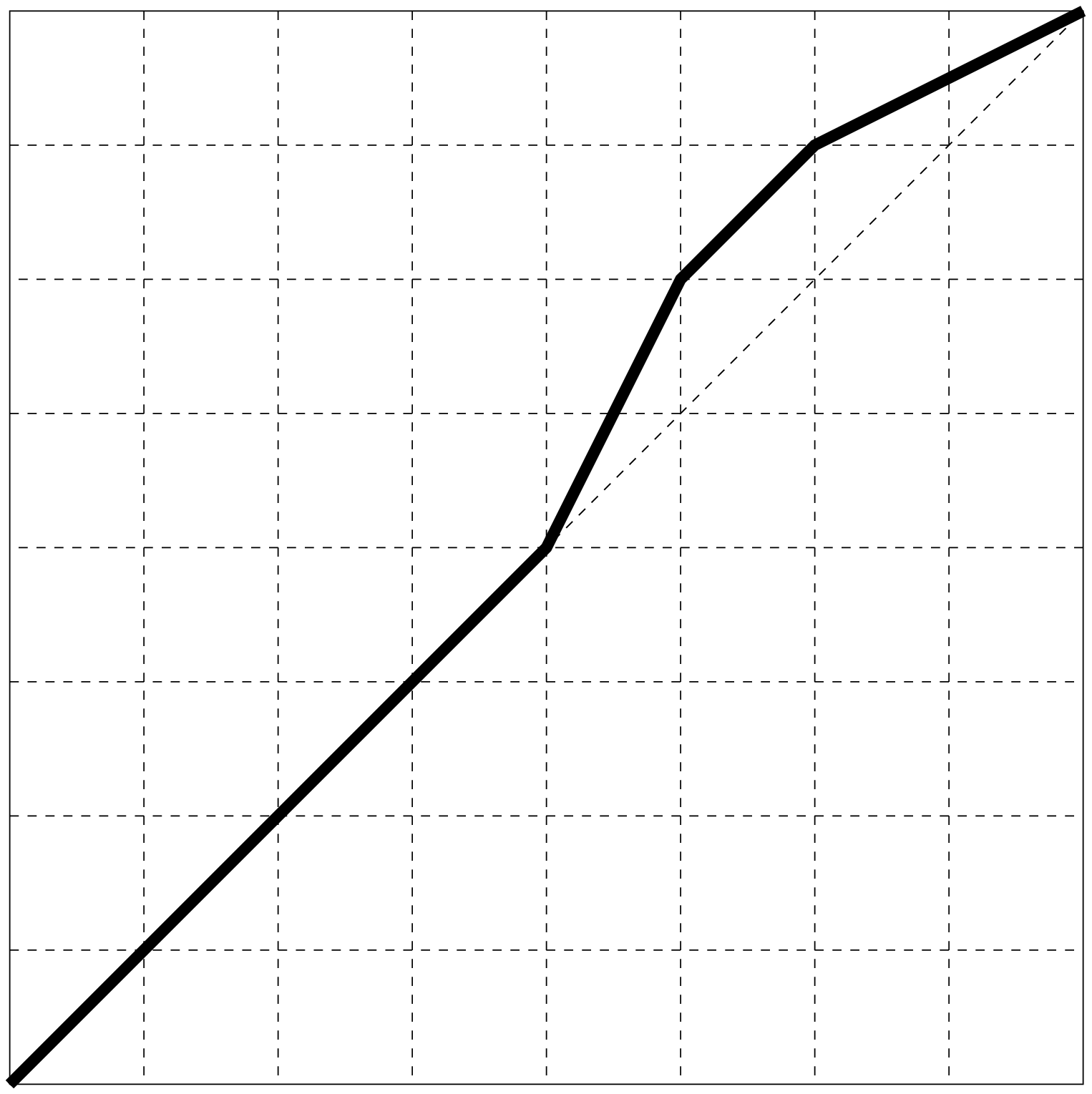}\\
$f_0$ \ \ \ \ \ \ \ \ \ \ \ \ \ \ \ \ \ \ \ \ \ \ \ \ \ \ \ \ \ \ \ \ \ \ \ \ \ \ \ \ \ \ \ \ $f_1$
\end{center}

Note that in the above, composition and evaluation of functions in $F$ will
be written in word order.  In other words, if $f,g \in F$ and $t \in
[0,1]$ , then, $tf = f(t)$, $fg = g \circ f$, and $f^g = g^{-1}fg = g
\circ f \circ g^{-1} $.  

One can check, using the convention above, that $f_0\sim x_0$ and
$f_1\sim x_1$ satisfy the relevant relations from the second
presentation.  It is well known that the second presentation is derived
from the first (see \cite{CFP} Theorem 3.4).  The fact that $f_0$ and
$f_1$ generate all of the claimed functions in $F$ (as a group of
homeomorphisms) is Corollary 2.6 in \cite{CFP}.
(Note that our functions $f_0$ and $f_1$ are the inverses of the
homeomorphisms they use.)

Given a homeomorphism $f:[0,1]\to [0,1]$, we will denote by $\supp(f)$
the {\em{support of $f$}}, wehre we take this set to be the set of all
points in $[0,1]$ which are moved by the action of $f$.  That is
\[
\supp(f)=\left\{x\in [0,1]|xf\neq f\right\}.
\]
(Note that this is different from the definition used in analysis,
where a closure is taken.)

The fact that elements of $F$ have piecewise-linear graphs that admit
only finitely many breaks in slope immediately implies that if $f \in
F$, then $\supp(f)$ is a finite union of disjoint open intervals.  We
will call each of these disjoint open intervals an {\em{orbital}} of
$f$.

\section{Previous Results}
Here we mention several lemmas necessary for our proof whose
statements and proofs are spread throughout the literature. (If we do
not give an indication of where a lemma may be found in the
literature, then the lemma is standard and simple, and its proof may
be taken as an exercise for the reader.)

\bl\label{disjointSupport}  

If $f$ and $g$ are set functions where the support of $f$ is disjoint
from the support of $g$, then $f$ and $g$ commute.

\el

\bl \label{unionSupports}

Let $g$, $f \in F$.  Let $H$ be the subgroup of $F$ that is generated
by $f$ and $g$ and define
\[
\supp(H) = \left\{x\in [0,1]|xh\ne x \textrm{ for some }h\in
H\right\}.
\]
Then, $\supp(H) = \supp(g) \cup \supp(f)$.

\el 

The first point in the following lemma is essentially standard from
the theory of permutation groups.  It is stated (basic fact (1.1.a)) in a general form in \cite{BSPLR}).  The second point is Remark
2.3 in \cite{bpgsc}.

\bl
\label{conjugateOrbitals} 

Let $f \in F$ and let $g\in \textrm{Homeo}([0,1])$ be any homeomorphism of
the unit interval. Further suppose that $(a_1, b_1), (a_2, b_2),
\dots, (a_n, b_n)$ are the orbitals of $f$.  Under these assumptions

\be
\item the orbitals of $f^g$ are exactly $(a_1g,
b_1g), (a_2g, b_2g), \dots, (a_ng, b_ng)$, and
\item if $g$ is orientation-preserving and piecewise-linear then for
every $i$, the derivative from the right of $f$ at $a_i$ equals the
derivative from the right of $f^g$ at $a_ig$ and the derivative from
the left of $f$ at $b_i$ equals the derivative from the left of $f^g$
at $b_ig$.  \ee \el

The first part of the following lemma is immediate from the
definitions, while the second part is essentially a restatement of Lemma
3.4 in \cite{BSPLR}.

\bl\label{transitiveOrbital} 

If $(a,b)$ is an orbital of $f \in F$ and if $c
\in (a,b)$, then 

\be
\item for all $m \in Z$, $cf^m \in (a,b)$ and

\item for any $\varepsilon > 0$, there is an $n \in Z$ so that both
$a < cf^n < a + \varepsilon$ and $b - \varepsilon < cf^{-n} < b$.  

\ee
\el

Given a group $G$ of orientation preserving homeomorphisms of $[0,1]$,
a set $X\subset [0,1]$, and a positive integer $k$, we say that
\emph{$G$ acts $k$-transitively over $X$} if given any two sets
$x_1<x_2<\ldots<x_k$ and $y_1<y_2<\ldots<y_k$ of points in $X$, there
is a $g\in G$ so that $x_ig = y_i$ for all indices $i$.

The following are restatements of Lemma 4.2 and Theorem 4.3 from
\cite{CFP}.  
\bl \label{kTransitive}

R. Thompson's group $F$ acts $k$-transitively over the diadic
rationals in $(0,1)$, for all positive integers $k$.

\el

\bl \label{npnaq}

$F$ has no proper non-abelian quotients.  

\el

In particular, if we can find an $f$ and $g$ where 
\[
[fg^{-1},g^{f}]=[fg^{-1}, g^{f^2}]=1
\]
 and $f$ and $g$ do not commute, then $f$ and $g$ generate a group
that is isomorphic to $F$.  We need two more standard facts about $F$
(this is a combination of Theorems 4.1 and 4.5 in \cite{CFP}).

\bl \label{FPrime} 

The group $F' = [F,F]$, the commutator subgroup of $F$, is simple.
Furthermore, $F'$ consists of all of the functions $f \in F$ such that
both $f'(0) = 1$ and $f'(1) = 1$.

\el

The final lemma is contained in the second author's thesis \cite{WDiss}.

\bl\label{BLemma}

If $G\leq F$ and $G \cong F$, then there are generators $g_0$, $g_1
\in G$ such that $\langle g_0,\, g_1 \rangle = G$, and for every
orbital $A$ of $g_0$, if $B$ is an orbital of $g_1$, then either $A
\cap B = \emptyset $ or $B \subseteq A$.

Furthermore, for the same functions $g_0$ and $g_1$ as described
above, if $A = (a_1, a_2)$ is an orbital of $g_0$ but not $g_1$ and
$A$ is not disjoint from the support of $g_1$, then there is an
$\varepsilon > 0$ such that either

\be

\item $g_0$ and $g_1$ are equal in
the interval $(a_1, \ a_1 + \varepsilon)$ and $(a_2 - \varepsilon,
a_2)$ is disjoint from the support of $g_1$,  or

\item $g_0$ and $g_1$ are equal in the interval $(a_2 - \varepsilon,
a_2)$ and $(a_1, a_1 + \varepsilon)$ is disjoint from the support
of $g_1$.  

\ee
\el

\section{Properties of the finite index subgroups of F}
Here we derive some nice properties of the finite index subgroups of
$F$.  In particular, we explore their relationships with $F'$, and we
examine the extent of their supports.

We begin with a simple lemma about infinite simple groups.
\bl \label{FIInfSimple}

Infinite simple groups do not admit proper subgroups of finite
index.

\el

\emph{Proof:} 

Let $G$ be an infinite simple group and let $H$ be a finite index
subgroup of $G$.  The right cosets $\{He, Hg_2, ..., Hg_n \}$ form a
set that $G$ acts on by multiplication on the right (here we are
denoting the identity of $G$ by $e$).

The action induces a homomorphism from $G$ to the symmetric group on
$n$ letters.  Since the codomain of this homomorphism is a finite
group, the kernel must be non-trivial.  Since $G$ is simple the kernel
must be all of $G$.  Now, if $n \ne 1$, then we can assume that $Hg_2
\ne He = H$.  But now $H = He=H\cdot(g_2g_2^{-1}) = (Hg_2)\cdot
g_2^{-1} = Hg_2$ (the last equality follows as the action is trivial).
Thus, $n=1$ and $G=H$.

\qquad$\diamond$

Here we have the first lemma from the introduction.
{\flushleft\bf Lemma \ref{FPrimeInFIF}}
{\it
If $H\leq F$ is a finite index subgroup of $F$, then $F'\leq H$.
}

{\em{Proof:}} 

Let $H$ be a finite index subgroup of $F$. The group $H \cap F'$ must be
finite index in $F'$, which is an infinite simple group by Lemma
\ref{FPrime}.  Now, by Lemma \ref{FIInfSimple}, $F' \subseteq
H$.

\qquad $\diamond$

We can now prove Lemma \ref{FIFNormal}.

{\flushleft\bf Lemma \ref{FIFNormal}} {\it Suppose $H$ is a finite
index subgroup of $F$, then $H\trianglelefteqslant F$.  
}

{\em{Proof:}} 

Suppose that $H$ is not normal in $F$.  Then there is an $f \in F$ so
that $f^{-1}Hf \ne H$.  In particular, there is an $h \in H$ so that
$f^{-1}hf \notin H$.  This last implies that $h^{-1}(f^{-1}hf) \notin
H$.  But $h^{-1}f^{-1}hf = [h^{-1},f^{-1}] \in F'$.  Since Lemma
\ref{FPrimeInFIF} assures us that $F'\subseteq H$, we have a
contradiction.

\qquad$\diamond$

Also, we are in a good position to prove the following.

{\flushleft{\bf Lemma \ref{innerOuter}}}

{\it

Suppose $H$ is a finite index subgroup of $F$.  Then there exists
a unique maximal rectangular subgroup $\inner(H)$ of $H$ and a unique
minimal rectangular group $\outer(H)$ containing $H$.  

}

{\it Proof:}

Let $H$ be a finite index subgroup $F$ and suppose $\fk{a}{b}\leq H$
and $\fk{c}{d}\leq H$.  Let $r = \gcd(a,c)$ and $s = \gcd(b,d)$.  We
can use a finite product of elements from $\fk{a}{b}$ anf $\fk{c}{d}$
to build an element $f$ with $\phi{f} = (r,0)$, and likewise, we can
build an element $g$ with $\phi(g) = (0,s)$.  Now, using Lemma
\ref{FPrimeInFIF} it is immediate that $\fk{r}{s}\leq H$.  In
particular, any finite index subgroup of $F$ has a unique, maximal
rectangular subgroup.

$F = \fk{1}{1}$ is a rectangular subgroup of $F$ which contains $H$,
and it is easy to see that the intersection of any two rectangular
subgroups of $F$ is again a rectangular subgroup of $H$, in
particular, the intersection of all of the rectangular subgroups of
$F$ which contain $H$ produces a unique minimal rectangular group
containing $H$.

\qquad$\diamond$

We now pass to some further useful lemmas not mentioned in the
introduction.

\bl \label{FIFSupports}

 If $H$ is finite index in $F$, then

\be

\item $\supp(H) = (0,1)$, and
\item there are $h_1$, $h_2 \in H$ so that $\supp(h_1) = (b,1)$ and
$\supp(h_2) = (0,a)$, for some $0<a \leq 1$ and $0 \leq b < 1$.

\ee
\el

{\em{Proof:}} 

(1) By the proof of Lemma \ref{FIFNormal}, $F'\leq H$,
and $\supp(F') = (0,1)$.

(2) Suppose that for all $h \in H$, $h'(1) = 1$.  Then, for all
$g_k \in Hf_0^{k}$, $(g_k)'(1) = \frac{1}{2^k}$.  In particular, we
have just found infinitely many distinct right cosets of $H$ in $F$.

A similar argument shows there is an $h \in H$ with $\supp(h) =
(0,a)$.

\qquad$\diamond$

\section{Finite Index Subgroups of F that are Isomorphic to F}
\label{FIFIsoF}
Consider the functions $g_0$ and $g_1$ specified by their sets of
breaks as follows:
\[
g_0 \textrm{ has breaks }
\left\{
\left(\frac{3}{8},\frac{3}{8}\right),\,
\left(\frac{1}{2},\frac{5}{8}\right),\,
\left(\frac{5}{8},\frac{3}{4}\right),\,
\left(\frac{7}{8},\frac{7}{8}\right)
\right\}
\]
\[
g_1 \textrm{ has breaks }
\left\{
\left(\frac{3}{8},\frac{3}{8}\right),\,
\left(\frac{7}{16},\frac{1}{2}\right),\,
\left(\frac{1}{2},\frac{9}{16}\right),\,
\left(\frac{5}{8},\frac{5}{8}\right)
\right\}
\]
These functions have graphs as below.

\begin{center}
\includegraphics[scale=0.25]{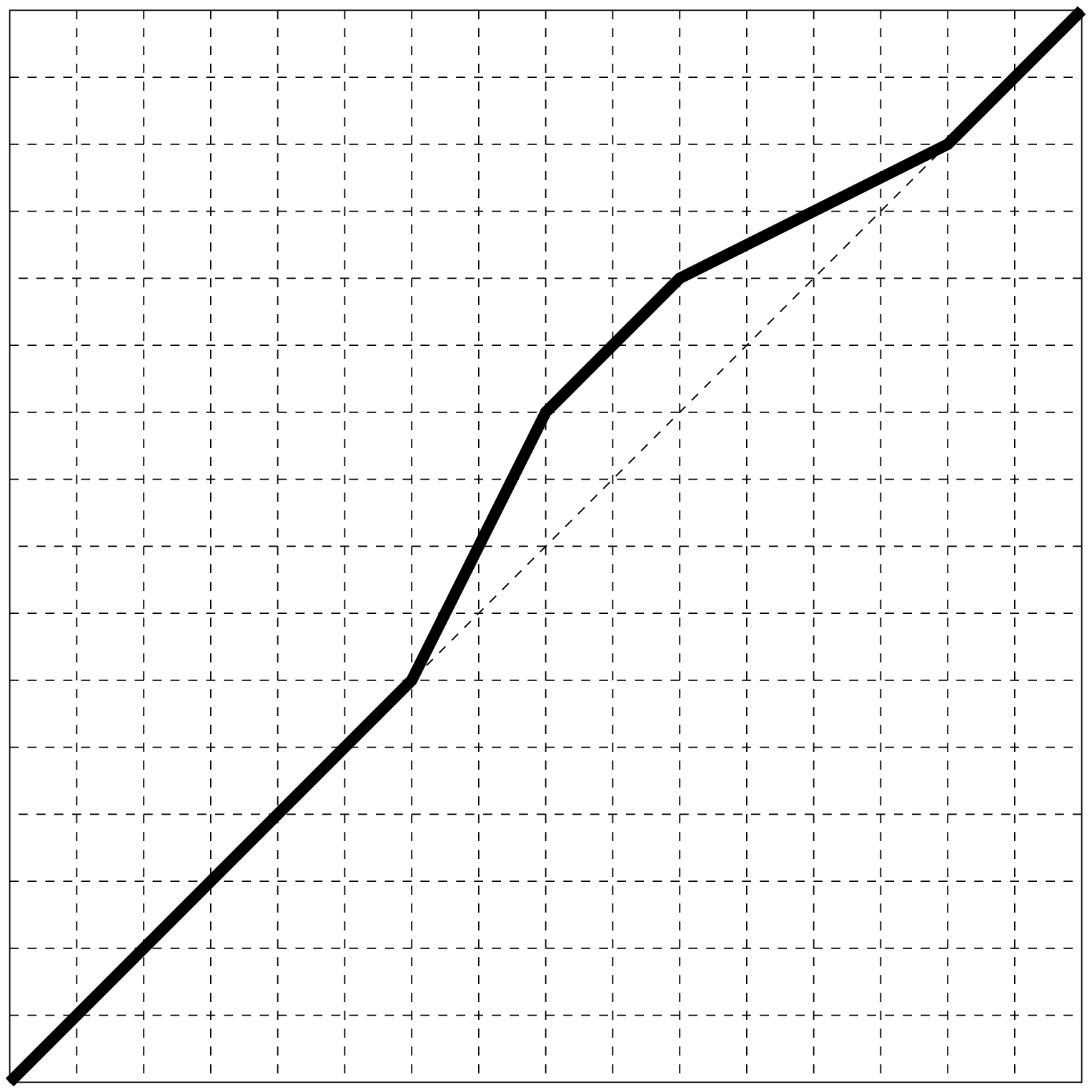} \ \ \ \ \ \ \ \ \ \ \ \ \ \ \ \ 
\includegraphics[scale=0.25]{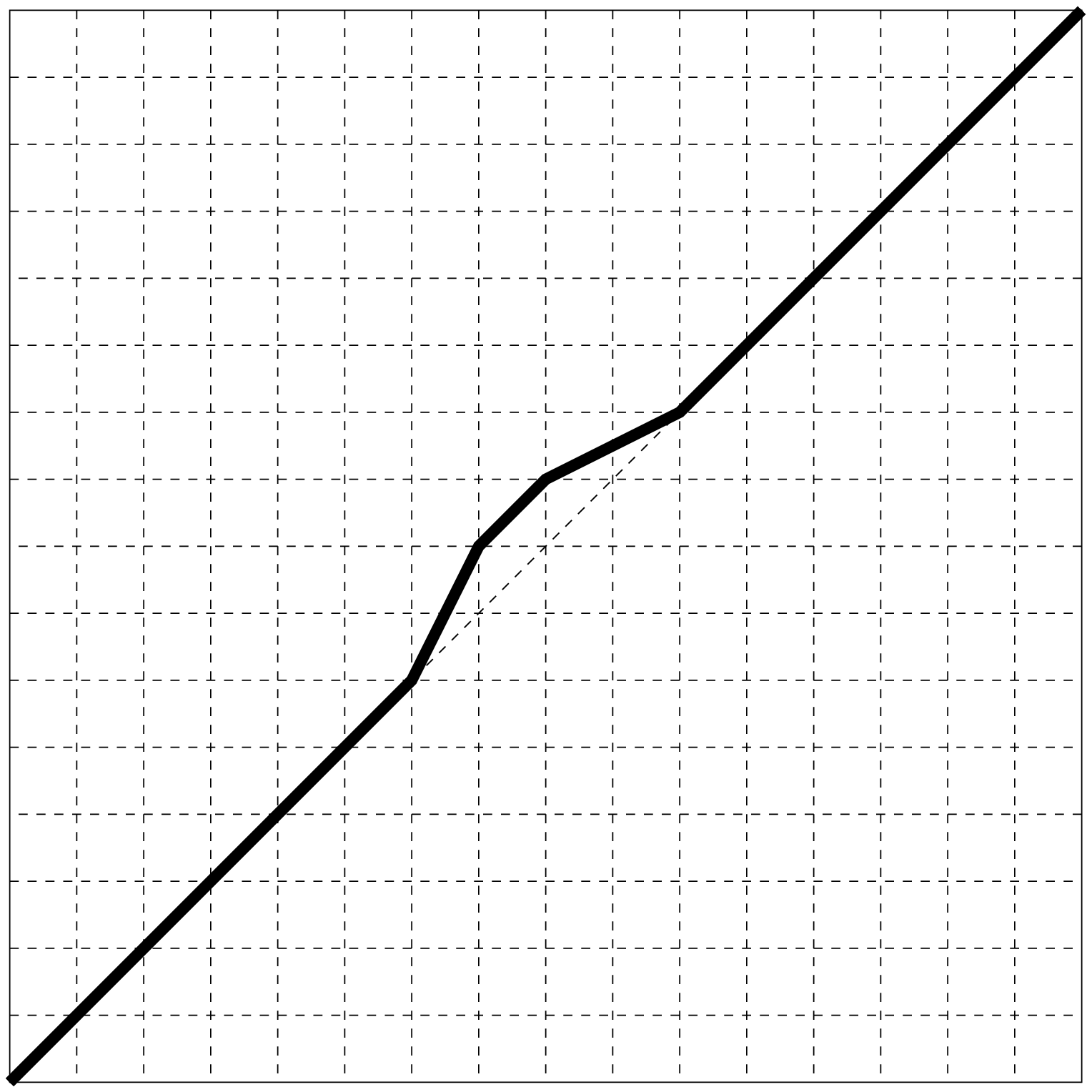}\\
$g_0$ \ \ \ \ \ \ \ \ \ \ \ \ \ \ \ \ \ \ \ \ \ \ \ \ \ \ \ \ \ \ \ \ \ \ \ \ \ \ \ \ \ \ \ \ $g_1$
\end{center}

\bl \label{fullBox}

  Let $g_0$ and $g_1$ be the functions in $F$ that are defined above.
Then $\langle g_0, g_1 \rangle$ 
\be
\item consists of every element of $F$ whose support is contained in
the interval $ [\frac{3}{8}, \frac{7}{8}]$, and
\item is isomorphic with $F$.  
\ee 
\el

{\em{Proof:}} 

We only need show the first point.  The second point will then follow
since by Lemma 4.4 in \cite{CFP} the subset of elements of $F$ with
support in $[a,b]$ where $a$ and $b$ are diadic rationals with $b-a$
an integral power of two is conjugate by a linear homeomorphism of $R$
to produce exactly $F$.

We explicitly build the linear conjugator of Lemma 4.4 in \cite{CFP}.

Consider the homeomorphism $\omega:R\to R$ defined by $t \mapsto
\left( \frac{8t-3 }{4} \right)$.  This homeomorphism sends
$[\frac{3}{8}, \frac{7}{8}]$ linearly to $[0,1]$, and it induces an
isomorphism $\psi: \langle g_0, g_1 \rangle \to H$ for some subgroup
$H\leq Homeo(R)$.  (Here, we are considering elements of $F$ to be
homeomorphisms from $R$ to $R$, by using the unique extension of any
element of $F$ by the identity map away from $[0,1]$).  The function
$\psi$ can be thought of as a restriction of the inner automorphism of
$Homeo(R)$ produced by conjugation by $\omega$.  

From here out, we will refer to $\langle g_0,g_1\rangle$ as $\Gamma$.

If we restrict $\psi$ less (potentially, depending on the size of
$\Gamma$), and take the preimage of $F$ under $\psi^{-1}$, then Lemma
4.4 in \cite{CFP} tells us that $\psi^{-1}(F) = \Upsilon \cong F$,
where $\Upsilon$ consists of all graphs of $F$ with support in $[3/8,
7/8]$.

Since $\omega$ is linear, we can understand $\psi$ by considering how
the map $\omega$ moves the breaks of any element in $\langle g_0,
g_1\rangle$. If $(p_i,q_i)$, $1 \leq i \leq n$, are the breaks of $g
\in \langle g_0, g_1 \rangle $, then $g \psi$ is the unique
piecewise-linear element of $Homeo(R)$ whose breaks are $( \frac{8p_i
-3 }{4}, \frac{8q_i -3 }{4} )$ which acts as the identity near $\pm
\infty$.

Now, one can check directly that $g_0 \psi = f_0$ and $g_1 \psi =
f_0^2f_1^{-1}f_0^{-1}$.  So $\langle g_0 \psi, g_1 \psi \rangle =
\langle f_0 , \ f_0^2f_1^{-1}f_0^{-1} \rangle = \langle f_0 , \
f_0^{-2} (f_0^2f_1^{-1}f_0^{-1}) f_0 \rangle = \langle f_0 , \
f_1^{-1} \rangle = \langle f_0 , \ f_1 \rangle = F $.  In particular,
$\psi(\Gamma) = F$, hence $\Upsilon = \Gamma$, and $\Gamma\cong F$.

\qquad$\diamond$

\begin{center}
\includegraphics[scale=0.25]{f0} \ \ \ \ \ \ \ \ \ \ \ \ \ \ \ \ 
\includegraphics[scale=0.25]{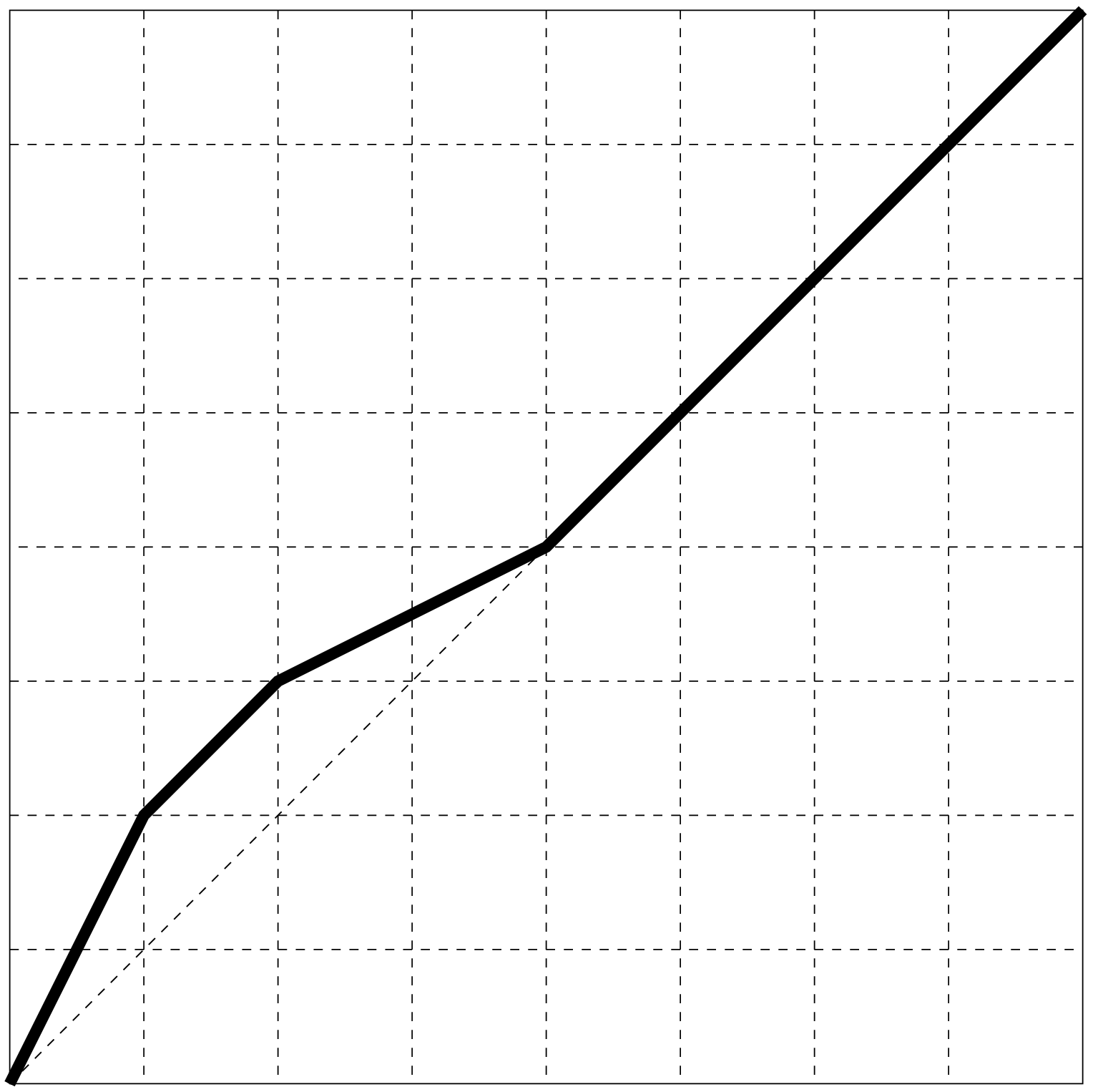}\\
$g_0 \psi$ \ \ \ \ \ \ \ \ \ \ \ \ \ \ \ \ \ \ \ \ \ \ \ \ \ \ \ \ \ \ \ \ \ \ \ \ \ \ \ \ \ \ $g_1 \psi$
\end{center}

We are now ready to prove the first of our main theorems.  For the
following, we need to recall the $\kab$ groups:

\[
K_{(a , b)} = \left\{  h \in F \  | \  \exists m,n \in Z \
  \mbox{s.t.} \  h'(0) = (2^{a})^n \ \mbox{and} \  h'(1) = (2^{b})^m
\right\}
\]

 where both $a$ and $b$ are non-zero integers.

\vspace{.1 in}

{\bf Theorem \ref{KabIsoF}}  {\it Let $H$ be a finite index subgroup of $F$.  $H$ is isomorphic to $F$
if and only if $H = K_{(a,b)}$ for some $a,b \in N$.
}

{\em{Proof:}} \ \ $( \Longleftarrow )$: Fix $a$ and $b$ in $N$.  We
will build generators $y_0$ and $y_1$ for for $ K_{(a,b)}$.  First, we
will define $y_0 \in K_{(a ,b )}$ over a finite collection of points
as follows:\\ If $a = 1$, then let $(a_1,b_1) = (\frac{1}{16} ,
\frac{1}{8})$.  If $a \not= 1$, then let $(a_1,b_1) =
(\frac{1}{2^{2a}} , \frac{1}{2^{a}})$.\\ Let $(a_2, b_2) =
(\frac{1}{8}, \frac{3}{8})$. Let $(a_3, b_3) = (\frac{5}{8},
\frac{7}{8})$.  If $b = 1$, then let $(a_4,b_4) = (\frac{7}{8} ,
\frac{15}{16})$.  If $b \not= 1$, then let $(a_4,b_4) =
(1-\frac{1}{2^{b}} , 1- \frac{1}{2^{2b}})$.

\underline{Filling in the definition of $y_0$}.

Extend the definition of $y_0$ by making it linear from $(0,0)$ to
$(a_1,b_1)$, affine and with slope one from $(a_2,b_2)$ to
$(a_3,b_3)$, and affine from $(a_4,b_4)$ to $(1,1)$.  All slopes
involved so far are integral powers of two, and the set os $a_i$'s and
$b_i$'s are all diadic rationals, so $y_0$ still has the potential to
be extended to an element of $F$.

We can now pick some diadic rational pairs $(c_1,d_1)$ and $(c_2,d_2)$
with $a_1<c_1<c_2<a_2$ and $b_1<d_1<d_2<b_2$ so that the ratios
\[
\frac{d_1-b_1}{c_1-a_1}\qquad\textrm{ and }\qquad\frac{b_2-d_2}{a_2-c_2}
\]
both produce integral powers of $2$ (not equal to the values of the
slope of $y_0$ near zero, or to the value $1$, the slope of $y_0$ over
$(a_2,a_3)$), and where the line segments from $(a_1,b_1)$ to
$(c_1,d_1)$ and from $(c_2,d_2)$ to $(a_2,b_2)$ (which we will be
adding to the definition of $y_0$) do not cross the line $y = x$.  We
can now extend the definition of $y_0$ from $0$ to $c_1$ and from
$c_2$ to $a_3$ so that over each interval, $y_0$ admits precisely one
breakpoint (over $a_1$ and $a_2$ respectively), and the graph of $y_0$
determined so far stays well above the line $y = x$.  By Lemma
\ref{kTransitive}, there is an element $\zeta$ of $F$ which sends the
list of points $(0,a_1,c_1,c_2,a_2,a_3,a_4,1)$ to the list
$(0,b_1,d_1,d_2,b_2,b_3,b_4,1)$.  Assume we have previously expanded
these lists as necessary with many diadic points imbetween $c_1$ and
$c_2$ and correspondingly many diadic points between $d_1$ and $d_2$,
(all new points roughly evenly spaced out) so that the graph of
$\zeta$ cannot intersect the line $y=x$.  We can now define $y_0$ over
the interval $(c_1,c_2)$ to agree with $\zeta$.  The element $y_0$ is
now defined over the intervals $(0,a_3)$ and $(a_4,1)$.

We can fill in the definition of $y_0$ with similar care over the
region $(a_3,a_4)$ (choose diadics $c_3$ and $c_4$ with
$a_3<c_3<c_4<a_4$ in a fashion similar to our choices of $c_1$ and
$c_2$, then connect over the region $(c_3,c_4)$ by some random
appropriate element of $F$ which does not touch the line $y = x$) to
finally get an element $y_0$ in $F$ which

\be
\item is linear over $(0,a_1)$, $(a_2,a_3)$, and $(a_4,1)$, and 
\item has breakpoints including $(a_1,b_1)$, $(a_2,b_2)$,
$(a_3,b_3)$, and $(a_4,b_4)$, and
\item does not intersect the line $y=x$.
\ee

Note that while $y_0$ is defined everywhere, it is not completely
determined over $(c_1,c_2)$, and it is not completely determined over
the similar interval $(c_3,c_4)$ in $(a_3,a_4)$ (although it is
roughly controlled in both locations).

\underline{Construct $y_1$ as follows}:
\[
ty_1 = \left\{
\begin{array}{ccc}
   t    & : & t \leq 3/8    \\
   2t- (3/8) & : &  3/8 \leq t \leq 5/8     \\
    ty_0    & : &  5/8 \leq t \leq 1    \\
\end{array} \right.
\]

{\textsc{sub-claim \ref{KabIsoF}.1:}} $\kab \trianglelefteq F$.

\vspace{8pt}

{\em{Proof of \ref{KabIsoF}.1:}} Let $g,h \in \kab$.  Suppose $g'(0) =
(2^a)^m$ and $h'(0) = (2^a)^n$.  Since all elements of $F$ are linear
in a neighborhood of $0$, then the chain rule for derivatives from the
right applies.  In particular $(gh)'(0) = (2^a)^{m+n}$.  Similarly,
$(gh)'(1) = (2^b)^{p+q}$, where $g'(1) = (2^b)^p$ and $h'(1) =
(2^b)^q$.  So $\kab$ is a subgroup of $F$.

Let $f \in F$.  From Lemma \ref{conjugateOrbitals}, it must be the case that $(g^f)'(0)
= (2^a)^m$ and $(g^f)'(1) = (2^b)^p$.  So $g^f \in \kab$.  Thus $\kab
\trianglelefteq F$.
 
\vspace{8pt}

{\textsc{sub-claim \ref{KabIsoF}.2:}} $\kab$ is a finite index subgroup of F.

\vspace{8pt}

{\em{Proof of \ref{KabIsoF}.2:}} Let $f \in F$.  The slope of $f$ near
$0$ is $2^p$ and the slopes of elements of $\kab$ near $0$ is $2^{a
n}$ for $n \in Z$.  Then the slopes of elements of $f \kab$ near $0$
is $2^{a n + p}$.  The division algorithm gives us that since $n \in
Z$, there are exactly $a$ different cosets of $K_{(a , 1)}$.
Similarly, there are exactly $b$ different cosets for $K_{(1 , b)}$.
Since $K_{(a , b)} = K_{(a , 1)} \cap K_{(1 , b)}$, then are at most
$ab$ distinct cosets for $\kab$ in $F$.

\vspace{8pt}

{\textsc{sub-claim \ref{KabIsoF}.3:}} $Y = \langle y_0, y_1 \rangle \cong F$. 

\vspace{8pt}

{\em{Proof of \ref{KabIsoF}.3:}} $y_0$ and $y_1$ have been constructed
specifically to have orbitals of certain products of these functions
to be disjoint. Since $y_0|_{[\frac{5}{8}, 1]} = y_1|_{[\frac{5}{8},
1]} $, and as both functions have graphs above the line $y=x$ in this
region, it must be the case that $\supp(y_0y_1^{-1}) = (0,
\frac{5}{8})$.  By Lemma \ref{conjugateOrbitals}, $\supp(y_1^{y_0}) =
(\frac{5}{8}, 1) $ and $\supp(y_1^{y_0^2}) = (\frac{7}{8}, 1)$.  By
Lemma \ref{disjointSupport}, $[y_0y_1^{-1},y_1^{y_0} ] = 1$ and
$[y_0y_1^{-1},y_1^{y_0^2} ] = 1$.  $y_0$ and $y_1$ do not commute
because $\frac{1}{4}y_0y_1y_0^{-1}y_1^{-1} =
\frac{1}{2}y_1y_0^{-1}y_1^{-1} = \frac{5}{8}y_0^{-1}y_1^{-1} =
\frac{3}{8}y_1^{-1} = \frac{3}{8} \not= \frac{1}{4}$.  So then by
Lemma \ref{npnaq}, $Y \cong F$.

\vspace{8pt}

{\textsc{sub-claim \ref{KabIsoF}.4:}} $Y = \langle y_0, y_1 \rangle = \kab$.

\vspace{8pt}

{\em{Proof of \ref{KabIsoF}.4}}: Note that $y_1'(0) = 1$, $y_1'(1) =
 y_0'(1)$,\\ $y_0'(0) = \left\{ \begin{array}{ccc} \frac{(1/8)}{
 (1/16)} = 2^1 & \mbox{if} & a = 1 \\ \frac{(1/2^a) }{(1/2^{2a}) } =
 2^a & \mbox{if} & a > 1\\
 \end{array}  \right.$ \  and \    $y_0'(1) = \left\{ \begin{array}{ccc}
 \frac{(1/16)}{ (1/8)} = 2^{-1}    & \mbox{if} & b = 1 \\
 \frac{(1/2^{2b}) }{(1/2^{b}) } = 2^{-b}   & \mbox{if} & b > 1\\
 \end{array}  \right.$.

So $y_0$ and $y_1$ are both in $\kab$ and $Y = \langle y_0, y_1
\rangle \subseteq \kab$.

We have carefully constructed $y_0$ and $y_1$ in such a way that even
though there are two intervals over which $y_0$ is not explicitly
known, the commutator function $[y_0, y_1]$ is completely determined.
Let us demonstrate this point.

Since $y_0|_{[5/8,1]} = y_1|_{[5/8,1]}$, then
$y_0y_1y_0^{-1}y_1^{-1}|_{[5/8,1]} = 1$.  Since $y_1|_{[0, 3/8]} = 1$
and $\frac{1}{8}y_0 = \frac{3}{8}$, then
$y_0y_1y_0^{-1}y_1^{-1}|_{[0,1/8]} = 1$.\\

The following line segments are taken linearly to each other.
\[
\left[ \frac{1}{8}, \frac{1}{4} \right] \stackrel{y_0}{\longmapsto} \left[ \frac{3}{8},
\frac{1}{2} \right] \stackrel{y_1}{\longmapsto} \left[ \frac{3}{8}, \frac{5}{8} \right]
\stackrel{y_0^{-1}}{\longmapsto} \left[ \frac{1}{8}, \frac{3}{8} \right]
\stackrel{y_1^{-1}}{\longmapsto} \left[ \frac{1}{8}, \frac{3}{8} \right] .
\]
\[
\left[ \frac{1}{4}, \frac{3}{8} \right] \stackrel{y_0}{\longmapsto}
  \left[ \frac{1}{2}, \frac{5}{8} \right] \stackrel{y_1}{\longmapsto}
  \left[ \frac{5}{8}, \frac{7}{8} \right]
  \stackrel{y_0^{-1}}{\longmapsto} \left[ \frac{3}{8}, \frac{5}{8}
  \right] \stackrel{y_1^{-1}}{\longmapsto} \left[ \frac{3}{8},
  \frac{1}{2} \right].
\]
\[
\left[ \frac{3}{8}, \frac{5}{8} \right] \stackrel{y_0}{\longmapsto}
  \left[ \frac{5}{8}, \frac{7}{8} \right]
\]

Since $y_1y_0^{-1}|_{[5/8,1]} = 1$, then $y_0y_1y_0^{-1}$ linearly
  maps $\left[ \frac{3}{8}, \frac{5}{8} \right] $ to $ \left[
  \frac{5}{8}, \frac{7}{8} \right]$, which is taken linearly by
  $y_1^{-1}$ to $[\frac{1}{2}, \frac{5}{8}]$.

 Now $y_0y_1y_0^{-1}y_1^{-1}$ contains the straight line segments from
 $(0,0)$ to $(\frac{1}{8}, \frac{1}{8} )$, from $(\frac{1}{8},
 \frac{1}{8} )$ to $(\frac{1}{4}, \frac{3}{8})$, from $(\frac{1}{4} ,
 \frac{3}{8})$ to $(\frac{3}{8}, \frac{1}{2})$, from $(\frac{3}{8},
 \frac{1}{2})$ to $(\frac{5}{8}, \frac{5}{8})$, and from
 $(\frac{5}{8}, \frac{5}{8})$ to $(1,1)$.

Since \ $\supp([y_0, y_1]) = (\frac{1}{8}, \frac{5}{8})$ and $y_0$ is
explicitely known in the interval $(\frac{1}{8}, \frac{5}{8})$, then we can
explicitely find $[y_0, y_1]^{y_0}$.  Also, since $\supp([y_0, y_1]^{y_0}) =
(\frac{3}{8}, \frac{7}{8})$ and $y_1^{-1}$ is explicitely known on
$(\frac{3}{8}, \frac{7}{8})$, then $[y_0, y_1]^{y_0y_1^{-1}}$ can also
be computed.   This computation gives that $[y_0, y_1]^{y_0} = g_0$ and
$[y_0, y_1]^{y_0y_1^{-1}} = g_1$, where $g_0$ and $g_1$ are the functions defined in the beginning of Section \ref{FIFIsoF}.  So then by Lemma \ref{fullBox}, $\langle g_0, g_1 \rangle$ contains every element
of $F$ that has support inside the interval $ (\frac{3}{8}, \frac{7}{8})$.

Since $g_0$ and $g_1$ are products of the funtions $y_0, y_1 \in \kab$, then
$Y$ and $\kab$ both contain every element of $F$ whose support is contained in
the interval  $ (\frac{3}{8}, \frac{7}{8})$.

Let $h \in F'$.  By Lemma \ref{FPrime}, there exists a $c,d \in (0,1)$ so that $\supp(h) \subseteq
(c,d)$.  By Lemma \ref{transitiveOrbital}, since $\supp(y_0) = (0,1)$, there is an $n \in Z$ so that
$\supp(h^{y_0^n}) \subseteq (\frac{3}{8}, dy_0^n)$, where $\frac{3}{8} < cy_0^n < dy_0^n
< 1$.  By Lemma \ref{transitiveOrbital}, since $\supp(y_1) = (\frac{3}{8}, 1)$, then there is an $m \in Z$ so that
$ \frac{3}{8} = \frac{3}{8}y_1^m < cy_0^ny_1^m < dy_0^ny_1^m < \frac{7}{8} $ and $\supp( (h^{y_0^n})^{y_1^m} ) \subseteq (\frac{3}{8}, \frac{7}{8}) $.  By the
previous argument, it must be the case that $h^{y_0^ny_1^m} \in Y$.  So
then $h = (h^{y_0^ny_1^m})^{y_1^{-m}y_0^{-n}} \in Y $.  So $F' \subseteq Y$.

Let $w \in \kab$.  There is an $n,m \in Z$ so that $w'(0) = 2^{an}$ and
$w'(1) = 2^{bm}$.  Since $w$, $y_0$, and $y_1$ are all linear functions in a neighborhoods of $0$ and $1$, then the chain
rule gives $(wy_0^{-n})'(0) = (2^{an})(2^a)^{-n} = 1$,  $(wy_0^{-n})'(1) =
(2^{bm})(2^{-b})^{-n} = 2^{b(m+n)}$,  $(wy_0^{-n}y_1^{m+n})'(0) = (1)(1)^{m+n}
= 1$, and   $(wy_0^{-n}y_1^{m+n})'(1) = 2^{b(m+n)}(2^{-b})^{m+n} = 1$.  So
$wy_0^{-n}y_1^{m+n} \in F' \subseteq Y \Rightarrow w \in Y$.  

Thus $\kab = Y \cong F$. 

\vspace{8pt}

 $( \Longrightarrow )$:  Assume that $H$ is a finite index subgroup of $F$ and
 $H \cong F$.

By Lemma \ref{FIFNormal}, $H \trianglelefteq F$.  By Lemma \ref{FPrimeInFIF}, $F'\leq H$ so $\supp(H) = (0,1)$.  There exists functions $h_0$ and $h_1$ so that $H = \langle h_0,
h_1 \rangle$ that satisfy the conditions listed in Lemma \ref{BLemma}.  One condition
in Lemma \ref{BLemma} is if $A$ is an orbital of $h_0$ and $B$ is an orbital of $h_1$,
then either $B \subseteq A$ or $B \cap A = \emptyset$.  This guarantees that if
$p$ is a fixed point of $h_0$, then $p$ is also a fixed point of $h_1$.  So
then the point $p$ will be a fixed point of the group $H$.  So $p \notin \
\mbox{supp}(H) = (0,1)$.  So either $p=0$ or $p=1$ and $\supp(h_0) = (0,1)$.

Since $h_0$ is not the identity near either $0$ or $1$, then there exist
non-zero integers $a$ and $b$ so that $h_0'(0) = 2^{a}$ and $h_0'(1) =
2^{b}$.  Lemma \ref{BLemma} also guarantees that either $h_1'(0) = 1$ and $h_1'(1) =
2^b$ or  $h_1'(0) = 2^a$ and $h_1'(1) = 1$.  Without loss of generality,
assume  $h_1'(0) = 1$ and $h_1'(1) = 2^b$.  We want to show that $H = \kab$.

$(\subseteq ):$  $h_0 \in \kab$ and $h_1 \in \kab$, so $H = \langle h_0, h_1
\rangle \subseteq \kab$. 

$(\supseteq ):$  Let $f \in \kab$.  So $f'(0) = 2^{an}$ and $f'(1) = 2^{bm}$
for some $n,m \in Z$.  Then, by the chain rule, $(fh_0^{-n})'(0) =
(2^{an})(2^a)^{-n}= 1 $ and   $(fh_0^{-n})'(1) = (2^{bm})(2^b)^{-n}=
2^{b(m-n)}$.  Also,  $(fh_0^{-n}h_1^{n-m})'(1) = 1(1)^{n-m}=1 $ and
$(fh_0^{-n}h_1^{n-m})'(1) = 2^{b(m-n)} (2^b)^{n-m} = 1$.  So then by Lemma
\ref{FPrime}, $fh_0^{-n}h_1^{n-m} \in F'$.  Since $H \trianglelefteq F$, then $F'
\subseteq H$.  So  $fh_0^{-n}h_1^{n-m} \in H$.  So then  $f =
fh_0^{-n}h_1^{n-m}h_1^{m-n}h_0^{n} \in H$.  So $H = \kab$.

\qquad$\diamond$

{\flushleft\bf Theorem \ref{FIsExtension}}\\
{\it Given any positive integers $a$ and $b$, $F$ can be regarded as a
non-split extension of $Z_a\times Z_b$ by $F$.  In particular,
there are maps $\iota$ and $\tau$ so that the following sequence is
exact.}

\[
\xymatrix{
1\ar[r]&F\ar[r]^{\iota}&F\ar[r]^{\!\!\!\!\!\!\!\!\!\!\!\!\tau} & 
Z_a\times Z_b\ar[r] &1.}
\]

\emph{Proof}:
This theorem is actually an immediate corollary to Theorem \ref{KabIsoF};
simply take $\iota$ to be the composition of the isomorphism from $F$
to $\kab$ with the inclusion map of $\kab$ into $F$.

\qquad$\diamond$

To prove Theorem \ref{FIFCorreFIZ2}, we will need to produce some analysis of the finite index subgroups of $Z^2$.

\section{Finite index subgroups of $Z^2$}
In this section we will prove two statements about the finite index
subgroups of $Z^2$.  While both of these statements could be taken as
straightforward exercises in an entry level graduate course in group
theory, we will include the proofs for completeness.

\bl\label{maximalKab}

Suppose $H$ is a finite index subgroup of $Z^2$.  Then there are
minimal positive integers $a$ and $b$ so that $\tilde{K}_{(a,b)}\leq
H$.  Further, if $\tk{c}{d}\leq H$ then $\tk{c}{d}\leq \tk{a}{b}$.

\el

\emph{Proof:}

$H$ is normal in $Z^2$ since $Z^2$ is abelian.  In particular, since
$H$ has finite index in $Z^2$, the group $T=Z^2/H$ is finite.
Therefore, there is a minimal positve integer $a$ so that $(a,0)\in
H$ and a minimal positive integer $b$ so that $(0,b)\in H$.  It is
now immediate that $\tk{a}{b}\leq H$.

Suppose $\tk{c}{d}\in H$.  Then $(c,0)\in H$.  The Euclidean
Algorithm now shows that $(j,0)\in H$, where $j = \gcd(a,c)$.  If
$a\nmid c$ we must have that $j< a$, which contradicts our choice of
$a$.  In particular, $a\mid c$ and $(c,0)\in \langle
(a,0)\rangle\leq \tk{a}{b}$.  A similar argument shows that
$(0,d)\in \tk{a}{b}$.  Since $\tk{c}{d}$ is generated by $(c,0)$
and $(0,d)$, we have that $\tk{c}{d}\leq \tk{a}{b}$.

\qquad$\diamond$

In the above lemma, we will call the group $\tk{a}{b}$ \emph{the maximal
$\tilde{K}$ group in $H$}.

\bl \label{cyclicQuotient}

Suppose $H$ is a finite index subgroup in $Z^2$ with maximal
$\tilde{K}$ group $\tk{a}{b}$.  The group $Q \cong H/\tk{a}{b}$ is
finite cyclic.

\el

\emph{Proof:} Thinking of $Z^2$ as a planar lattice, the points in $H$
not in $\tk{a}{b}$ are the points which will survive under modding $H$
out by $\tk{a}{b}$ to become non-trivial elements of $Q$.  Thus, we
can find $Q$ as a subgroup of points in the finite rectangular lattice
$L = Z_a\times Z_b$.  Furthermore, as $a$ and $b$ are minimal positive
so that $(a,0)\in H$ and $(0,b)\in H$, we must have that the only
intersection $Q$ will have with the vertical axis in $L$ (the points
of the form $(0,r)$) or with the horizontal axis in $L$ (the points of
the form $(r,0)$) is at the point $(0,0)$.

In particular, suppose $(r,s)$ and $(t,u)$ are points in $Q$.  If $j
\equiv \gcd(r,t)$, then we can again exploit the Euclidean Algorithm
to find integers $p$ and $q$ so that $p(r,s)+q(t,u) =(j,m)$ so that
$j$ divides both $r$ and $t$.  In $Z_a\times Z_b$ the point $(j,m)\in
Z^2$ becomes $(j,m_b)$.  Now there are positive integers $x$ and $y$
so that $x(j,m_b) = (r,xm_b)$ and $y(j,m_b) = (t,ym_b)$.  If
$xm_b\not\equiv s\mod b$ then $Q$ has an intersection with the
vertical axis in $Z_a\times Z_b$ away from $(0,0)$ and if
$ym_b\not\equiv u\mod b$ then $Q$ has an intersection with the
vertical axis of $Z_a\times Z_b$ away from $(0,0)$.  Since neither of
these intersections can exist, by the definitions of $a$ and $b$, we
see that $(r,s)$ and $(t,u)\in \langle (j,m_b)\rangle$ in $Q$.  In
particular, after a finite induction we see that $Q$ is cyclic.

\qquad$\diamond$

\section{The structure of the extension}
We have now done enough work so that Theorem \ref{FIFCorreFIZ2}
is transparent.  

{\flushleft\bf Theorem \ref{FIFCorreFIZ2}}
{\it \be
\item The map $\phi$ induces a one-one correspondence between the
finite index subgroups of $F$ and the finite index subgroups of $Z^2$.
\item Let $H$ be a finite index subgroup $H$ of $F$, with image
$\tilde{H}=\phi(H) \leq Z\times Z$.  There exist smallest positive
integers $a$ and $b$ with $\kab\leq H$.  Furthermore, if
$Q=\tilde{H}\slash\tilde{K}_{(a,b)}$, then $Q$ is finite cyclic, and
there are maps $\iota$, $\rho$, $\tilde{\iota}$ and $\tilde{\rho}$ so
that the diagram below commutes with the two rows being exact:
\ee
}
\[
\xymatrix { 
&  F\ar[d]_{\cong}  &  &  &  \\ 
1\ar[r]  &  K\ar[r]^{\iota} \ar[d]_{\phi|_K}  & 
H\ar[r]^{\rho}\ar[d]_{\phi|_H} & Q\ar@{=}[d]\ar[r] & 1\\ 
1\ar[r]  &  {\tilde{K}}\ar[r]^{\tilde{\iota}}  & 
{\tilde{H}}\ar[r]^{\tilde{\rho}}  &  Q\ar[r]  &  1.
}
\]

\vspace{.1 in}

{\it Proof:}
The first point follows from Lemma \ref{FPrimeInFIF} and the fact that the kernel of $\phi$ is $F'$.

The second point follows from a conglomeration of lemmas.  

The existence of minimal positive integers $a$ and $b$ (so that $\kab$
is maximal in $H$) follows from the existence of a maximal $\tk{a}{b}$
in $\tilde{H}$, which is lemma \ref{maximalKab}.

The fact that $Q$ is finite cyclic comes from Lemma \ref{cyclicQuotient}.

The isomorphism from $F$ to $K=\kab$ comes from Theorem \ref{KabIsoF}.

The map $\iota$ is the inclusion map of $\kab$ into $H$.  The map
$\tilde{\iota}$ is induced from the projection $\phi$. the map
$\tilde{\rho}$ is the natural quotient onto $Q$ of the image of
$\tilde{\iota}$ in $\tilde{H}$.  The bottom row is thus exact.  $\rho$
is the composition of the natural quotient of $H$ by the image of
$\iota$ followed by the isomorphism from $H/\iota(K)$ to
$\tilde{H}/\tilde{\iota}(\tilde{K}) = Q$, thus, the top row is exact,
and the diagram commutes.

\qquad$\diamond$

To prove Lemma \ref{characteristicKab} we will make use of Rubin's
Theorem.  The version we will quote is Theorem 2 in Brin's paper
\cite{BrinHigherV}.  That version is itself derived from Theorem 3.1
in the paper \cite{RubinWError} of Rubin, where in the statement of the
theorem, a technical hypothesis is inadvertently missing (see the
discussion of this in \cite{BrinHigherV}).

In order to state Rubin's Theorem, we will need to define some
terminology.  In this, we generalize the language of the definition of
\emph{locally dense} given in Brin's \cite{BrinHigherV}.  Our
generalization will have no impact on the content of our statement of
Rubin's theorem.  Suppose $X$ is a topological space and $H(X)$ is its
full group of homeomorphisms.  Suppose further that $K\leq H(X)$.
Given $W\subset X$, we will say \emph{$K$ acts locally densely over
$W$} if for every $w\in W$ and every open $U\subset W$ with $w\in U$,
the closure of

\[
\left\{ w\kappa |\kappa\in K, \kappa|_{(W-U)}=1_{(W-U)}\right\}
\]
contains some open set in $W$.  In particular, for each open $U$ in $Z$, the
subgroup of elements fixed away from $U$ has every orbit in $U$ dense
in some open set of $W$ in $U$.

We are now ready to state Rubin's theorem.  We give essentially the
statement given in \cite{BrinHigherV}, although we recast it in the
language of right actions.

\bt[Rubin] Let $X$ and $Y$ be locally compact, Hausdorff topological
spaces without isolated points, Let $H(X)$ and $H(Y)$ be the self
homeomorphism groups of $X$ and $Y$, respectively, and let $G\subseteq
H(X)$ and $H\subseteq H(Y)$ be subgroups.  If $G$ and $H$ are
isomorphic and both act locally densely over $X$ and $Y$,
respectively, then for each isomorphism $\varphi:G\to H$ there is a
unique homeomorphism $\gamma:X\to Y$ so that for each $g\in G$, we have
$g\varphi=\gamma^{-1} g\gamma$.  \et

In our case, and to apply Rubin's theorem to $F$ or subgroups of $F$,
we need to consider these groups to be groups of homeomorphisms of
$(0,1)$, instead of $[0,1]$. This comes from the simple fact that $F$
does not move $0$ or $1$ to produce a dense image in any open set!

Having made that (temporary) change to our definition of $F$ and its
subgroups, we are ready to apply Rubin's theorem to any such subgroup,
as long as it is locally dense in its action on $(0,1)$.  

In the discussion which follows, given $X\subset R$, we will use the
notation $D_X$ to denote the set $Z[1/2]\cap X$ of diadic rationals in
$X$.

\bl \label{FIFLocallyDense}
Finite index subgroups of $F$ act locally densely on $(0,1)$.
\el

{\it Proof:}

Suppose $H$ is finite index in $F$, and $x\in (0,1)$ and $U$ an open
neighborhood of $x$ in $(0,1)$.  Let $d_1$ and $d_2$ be two diadic
rationals in $U$ with $d_1<x<d_2$.  let $K$ be the subgroup of $F$
consisting of all the elements of $F$ with support in $(d_1,d_2)$. Let
$\alpha:R \to R$ be any piecewise-linear homeomorphism which is the
identity near $\pm \infty$ and which has all slopes integral powers of
$2$, and with all breaks occuring over the diadic rationals, and that
maps $d_1$ to $0$ and $d_2$ to $1$.  It is easy to build such a map,
and the reader may check that the inner automorphism of Homeo($R$)
generated by conjugation by $\alpha$ will take $K$ isomorphically to
$F$.

Now by an induction argument (for instance, as carried out in the
first paragraph of Section \S1. in \cite{CFP}), it is easy to see that
$\alpha$ takes $D_{(d_1,d_2)}$ to $D_{(0,1)}$ in an order preserving
fashion.  In particular, as $F$ is $k$-transitive on $D_{(0,1)}$ for any
positive integer $k$ (recall Lemma \ref{kTransitive}), we see that $K$
is $k$-transitive on $D_{(d_1,d_2)}$ for any positive integer $k$.

Now, if $x$ is a diadic rational, then the orbit of $x$ under $K$ is
dense in $(d_1,d_2)$, as $K$ acts transitively over
$D_{(d_1,d_2)}$, and $D_{(d_1,d_2)}$ is dense in $(d_1,d_2)$.

If $x$ is not diadic rational, then given any $\epsilon>0$, and any
$y$ in $(d_1,d_2)$, we can find four diadic rationals $x_1$,$x_2$,
$y_1$, $y_2\in (d_1,d_2)$ so that $x_1<x<x_2$ and $y_1<y<y_2$, and
where the $y_i$ are chosen epsilon-close to $y$.  Now there is some
element $\kappa$ in $K$ which throws $x_1$ to $y_1$ and $x_2$ to $y_2$
(since $K$ is $2$-transitive over $D_{(d_1,d_2)}$).  In particular,
$|y-x\kappa|<\epsilon$.  Hence, the orbit of $x$ is dense in
$(d_1,d_2)$.

Now, as $K\leq F'\leq H$, $H$ is locally dense over $(0,1)$.
\qquad$\diamond$

We will use Rubin's theorem to prove the final lemma from the introduction.

{\flushleft {\bf Lemma \ref{characteristicKab}}} 

{\it
Suppose $H$, $H'\in FIF$, $K=\inner(H)$,
$K' = \inner(H')$, and $\xi:H\to H'$ is an isomorphism. Then 

\be
\item $\xi(K)=K'$
\item $K$ is characteristic in $H$ and $K'$ is characteristic in $H'$, and
\item $\res(H) = \res(H')$.
\ee
}

{\it Proof:}

First, let us suppose $\vartheta:H\to H'$ is an isomorphism.

By Lemma \ref{FIFLocallyDense}, $H$ and $H'$ both act locally densely
on $(0,1)$.  In particular, Rubin's theorem tells us that there is a
homeomorphism $\gamma: (0,1)\to(0,1)$ so that for any $h\in H$,
$\vartheta(h) = \gamma^{-1}h\gamma \in H'$.

Now by Lemma \ref{conjugateOrbitals}, we see that the collection of
orbitals of $h'=\vartheta(h)$ is in bijective correspondence with the
orbitals of $h$. 

 Further, if $\gamma$ is orientation-preserving, any
orbital of $h$ which has end $e\in \left\{0,1\right\}$ becomes (under
the action of $\gamma$) an orbital of $h'$ with end $e$.  If $\gamma$ is
orientation-reversing, then any orbital of $h$ with end
$e\in\left\{0,1\right\}$ becomes an orbital of $h'$ with end $f\ne e$,
where $f\in\left\{0,1\right\}$.

Since $\vartheta$ is a homomorphism, a consequence of the above
paragraph is if $K = \fk{a}{b}$ for some positive integers $a$ and
$b$, then there are positive integers $c$ and $d$ with $\vartheta(K) =
\fk{c}{d}\leq K'\leq H'$.

The correspondence theorem now tells us that the maximal rectangular
groups of $H$ and $H'$ are mapped precisely to each other by
$\vartheta$, and we have point (1).

The second two points follow immediately.

Note that this argument provides a second proof that amongst the
finite index subgroups of $F$, only the rectangular groups are
actually isomorphic to $F$.  

\qquad $\diamond$

The lemma above provides the key ingredients for the proof of our
final theorem.

Recall the isomorphisms $\tau_{(a,b)}:\kab\to F$ from the introduction
(elements of $\kab$ with slope $(2^a)^s$ near zero are taken to elements of
$F$ with slope $2^s$ near zero, and elements of $F$ with slope $(2^b)^t$
near one are taken to elements of $F$ with slope $2^t$ near one), and
the map Outer which, given a finite index subgroup $H$ of $F$,
produces the smallest rectangular subgroup of $F$ that contains $H$.
With these maps in mind, and with the above lemma in hand, we are
finally ready to prove our last theorem.

{\flushleft {\bf Theorem \ref{isoClassification}:}}

{\it

Suppose $H$, $H'$ are finite index subgroups of $F$.  Let $a$, $b$,
$c$, $d$ be positive integers so that $\kab= \outer(H)$,
$K_{(c,d)}=\outer(H')$. $H$ is isomorphic with $H'$ if and only if
$\tau_{(a,b)}(H) = \tau_{(c,d)}(H')$ or $\tau_{(a,b)}(H) =
\rev(\tau_{(c,d)}(H'))$.  

}

{\it Proof:}

Suppose that $\vartheta:H\to H'$ is an isomorphism.  Lemma
\ref{characteristicKab} assures us that there is a well defined
positive integer $n$ so that $\res(H) = \res(H')=n$, and $\vartheta(K)
= K'$.  Let us further suppose that $\fk{r}{s}=\inner(H)$ and
$\fk{t}{u}=\inner(H')$.

Let $\tilde{H} = \phi(H)$, and $\tilde{H}' = \phi(H')$.  Consider the
translations of $Z^2$ generated by $(r,0)$ and $(0,s)$.  Since
$\tilde{H}$ is a group, the sets $\tilde{H}$, $\tilde{H}+(r,0)$, and
$\tilde{H}+(0,s)$ are the same.  In particular, we can consider the
image in the lattice $Z^2$ of $\tilde{H}$, restricting our view to the
rectangle $R$ of points with integer coordinates where the horizontal
coordinates range from $0$ to $r-1$ and the vertical coordinates range
from $0$ to $s-1$, and understand everything about the group
$\tilde{H}$.  $\tk{r}{s}$ only intersects $R$ at $(0,0)$, while there
are $n$ total intersections of $\tilde{H}$ with $R$, all obtained by translating different powers of some particular vector $(p,q)$ into $R$ (using $(r,0)$ and $(0,s)$).  Let $j = \mbox{gcd}(p,r)$.  So, the lowest column number that the image of the translated powers of $(p,q)$ in $R$ will appear in is column $j$.  Since $\tilde{H}$ intersects $R$ exactly $n$ times, it must be that case that $nj = r$ and the images of the translated powers of $(p,q)$ in $R$ will occur in columns $0$, $r/n$, $2r/n$, ... ,$(n-1)r/n$.  Similarly, $\mbox{gcd}(q,s) = s/n$ and the images of the translated powers of $(p,q)$ in $R$ will occur in
rows $0$, $s/n$, $2s/n$, ... ,$(n-1)s/n$.  Now, as $\tk{a}{b}$ is the smallest
rectangular group to contain $\tilde{H}$, we see that $a = r/n$ and $b
= s/n$.  A similar discussion shows that $c = t/n$ and $d = u/n$.
Stated another way, we have

\[
\frac{r}{a}=\frac{s}{b}=\frac{t}{c}=\frac{u}{d}=n.
\]

Now, consider the image of $H$ and $H'$ under the respective maps
$\tau_{(a,b)}$ and $\tau_{(c,d)}$.  The subgroups
$\fk{r}{s}=\inner(H)$ and $\fk{t}{u}=\inner(H')$ are both taken to
$\fk{n}{n}$. We will now assume that this is how $H$ and $H'$ started
out, and do all remaining work in these scaled versions of $H$ and
$H'$.

The isomorphism $\vartheta$ which is carrying $H$ to $H'$ must now
preserve the maximal rectangular subgroup $\fk{n}{n}$, by Lemma
\ref{characteristicKab}.  By Lemma \ref{FIFLocallyDense}, both $H$ and
$H'$ act locally densely on $(0,1)$, so by Rubin's theorem there is a
homeomorphism $\gamma$ so that for any $h\in H$, $\vartheta(h) =
\gamma^{-1}h\gamma\in H'$.  

Note that as $\gamma$ need not be piecewise-linear, we should be
concerned that conjugating by $\gamma$ might change slopes, as well as
potentially swapping coordinates.  

By Lemma \ref{characteristicKab} we know that $\fk{n}{n}=
\inner(H)=\inner(H')$ is being brought isomorphically to itself by
$\vartheta$.  Suppose $h\in H$ has an orbital $A$.  Denote by $E_A$
the set of ends of $A$ which are in the set $\left\{0,1\right\}$.  Now
consider $h' = \vartheta(h)$.  The element $h'$ has an orbital $B =
\gamma(A)$ by point (1) of Lemma \ref{conjugateOrbitals}.  Denote by
$E_B$ the ends of $B$ that are actually in the set
$\left\{0,1\right\}$.  Then as $\gamma$ preserves the set
$\left\{0,1\right\}$ we see that the cardinalities of $E_A$ and $E_B$
must be the same.

Now, by the result of the previous paragraph, and using the fact that
the $\vartheta$ takes $\fk{n}{n}$ isomorphically to itself, we see
that if $\gamma$ is orientation-preserving, we must have that $\gamma$
will send $(n,0)$ to $(n,0)$ and $(0,n)$ to $(0,n)$ in the induced map
from $\phi(H)\to \phi(H')$ (note that $(n,0)$ will not be taken to
$(-n,0)$, as conjugation by an orientation-preserving $\gamma$ will
preserve the local directions that points near zero and one move under
the action of $h$).  Similarly, if $\gamma$ is orientation-reversing,
then the reader can check that the action of $\gamma$ will send
$(n,0)$ to $(0, n)$ and $(0,n)$ to $(n,0)$, again considering the
induced map from $\phi(H)\to\phi(H')$.

 If $\gamma$ is orientation-reversing, then replace $H'$ by the
isomorphic copy $Rev(H')$, so that from here out we only need to argue
the case where our isomorphism $\vartheta$ appears to be the identity
after passing through the quotient map $\phi$.

Suppose $f\in H$.  $\phi(f) = (v,w)$, if and only if
$\phi(\vartheta(f))=(v,w)$.  But now as $H$ and $H'$ both contain the
commutator subgroup $F'$, and as they each contain an element which
has slopes $2^v$ and $2^w$ near zero and one respectively, we see that
both $H$ and $H'$ contain all of the elements of $F$ with $\phi(k) =
(v,w)$.  It is now immediate that $H = H'$.

Now let us suppose that $H$ and $H'$ are finite index subgroups of
$F$, and that $\fk{a}{b}= \outer(H)$ and $\fk{c}{d} = \outer(H')$.
Let us further suppose that the scaling maps $\tau_{(a,b)}$ and
$\tau_{(c,d)}$ have the property that $\tau{(a,b)}(H) =
\tau{(c,d)}(H')$ or $\tau_{(a,b)}(H)=\rev(\tau{(c,d)}(H'))$.  Since
the $\tau_{(*,*)}$ maps are isomorphisms, and the map $\rev$ is an
isomophism, we immediately see that $H$ and $H'$ are isomorphic.
\qquad$\diamond$

\section{Some examples}
In this section, we give some examples of finite index subgroups of
$F$, and consider them from the perspective of this paper.

  {\textsc{Example 1:}}\\ Let $H = \{ f \in F \ |
\ f'(0) = 2^{3n+5m} \ \mbox{and} \ f'(1) = 2^{7n+11m} \ \mbox{for
some} \ m,n \in Z\}$.\\ $H$ is a finite index subgroup of $F$ but $H$
is not isomorphic to $F$.

Let $\tilde{H} = \phi(H)$. Since $3, 5, 7, $ and $11$ are all odd, the
only possible elements of $\tilde{H}$ are $(\mbox{even}, \
\mbox{even})$ and $(\mbox{odd}, \ \mbox{odd})$.  If $n=5$ and $m=-3$,
then $(15-15, 35-33) = (0,2) \in \tilde{H}$.  If $n=-11$ and $m=7$,
then $(-33+35, -77+77) = (2,0) \in \tilde{H}$.  So then every
$(\mbox{even}, \ \mbox{even})$ is in $\tilde{H}$.  If $n=-3$ and
$m=2$, then $(1,1) \in \tilde{H}$.  Since $(1,1)$ and all
$(\mbox{even}, \ \mbox{even})$ are in $\tilde{H}$, then all
$(\mbox{odd}, \ \mbox{odd})$ are also in $\tilde{H}$.  So $\tilde{H}$
is index 2 in $Z \oplus Z$ and $H$ is index 2 in $F$.

To show that $H$ is not isomorphic to $F$, it is enough to show that $H \ne
\kab$ for any non-zero integers $a$ and $b$.

Assume that for some non-zero integers $a$ and $b$, $H = \kab$.  If $h
\in H$, there are integers $p$ and $q$ such that $\phi(h) = (ap, bq)$.

There is an $f \in H$ so that $\phi(f) = (3,7)$.  There is a $g \in H$
so that $\phi (g) = (5,11)$.  Now, there must be integers $p_1$ and
$p_2$ so that $ap_1=3$ and $ap_2=5$.  Thus $a = 1$.  Also, there must
be integers $q_1$ and $q_2$ so that $bq_1=7$ and $bq_2=11$. So $b=1$.
But $K_{(1,1)} = F \ne H$, so $H$ can not be isomorphic to $F$.

\qquad$\diamond$

In the above example, note that the maximal $\tk{a}{b}$ group was
$\tk{2}{2}$, which was proper in $\tilde{H}$, and that the quotient
$\tilde{H}/\tk{2}{2}\cong Z_2$.  In particular, $H$ is isomorphic to a
non-split extension of $F$ by $Z_2$, where the structure of the
extension is described by the structure of $\tilde{H}$ as an extension
of $\tk{2}{2}$ by $Z_2$.

  {\textsc{Example 2:}}\\ Let $l$,$r$,$f$, $f'\in F$, so that $\phi(l)
= (15,0)$, $\phi(r) = (0,15)$, $\phi(f) = (3,3)$, and $\phi(f') =
(3,6)$.  Suppose that $H = \langle F',l,r,f\rangle$ while $H'=\langle
F',l,r,f'\rangle$.

We see immediately that the maximal rectangular subgroups of $H$ and
$H'$ are $\fk{15}{15}$.  The minimal rectangular subgroups of $F$
containing $H$ and $H'$ are the same, namely $\fk{3}{3}$.  The
residues of $H$ and $H'$ are both $5$, but $H$ and $H'$ are not
isomorphic, as $\tau_{(3,3)}(H)\ne \tau_{(3,3)}(H')$ and
$\tau_{(3,3)}(H)\ne\rev(\tau_{(3,3)}(H'))$.

Below is included a diagram of the rectangle $R$ in $Z^2$ which demonstrates this non-equality.

\begin{center}
\includegraphics[height = 140 pt, width = 200 pt]{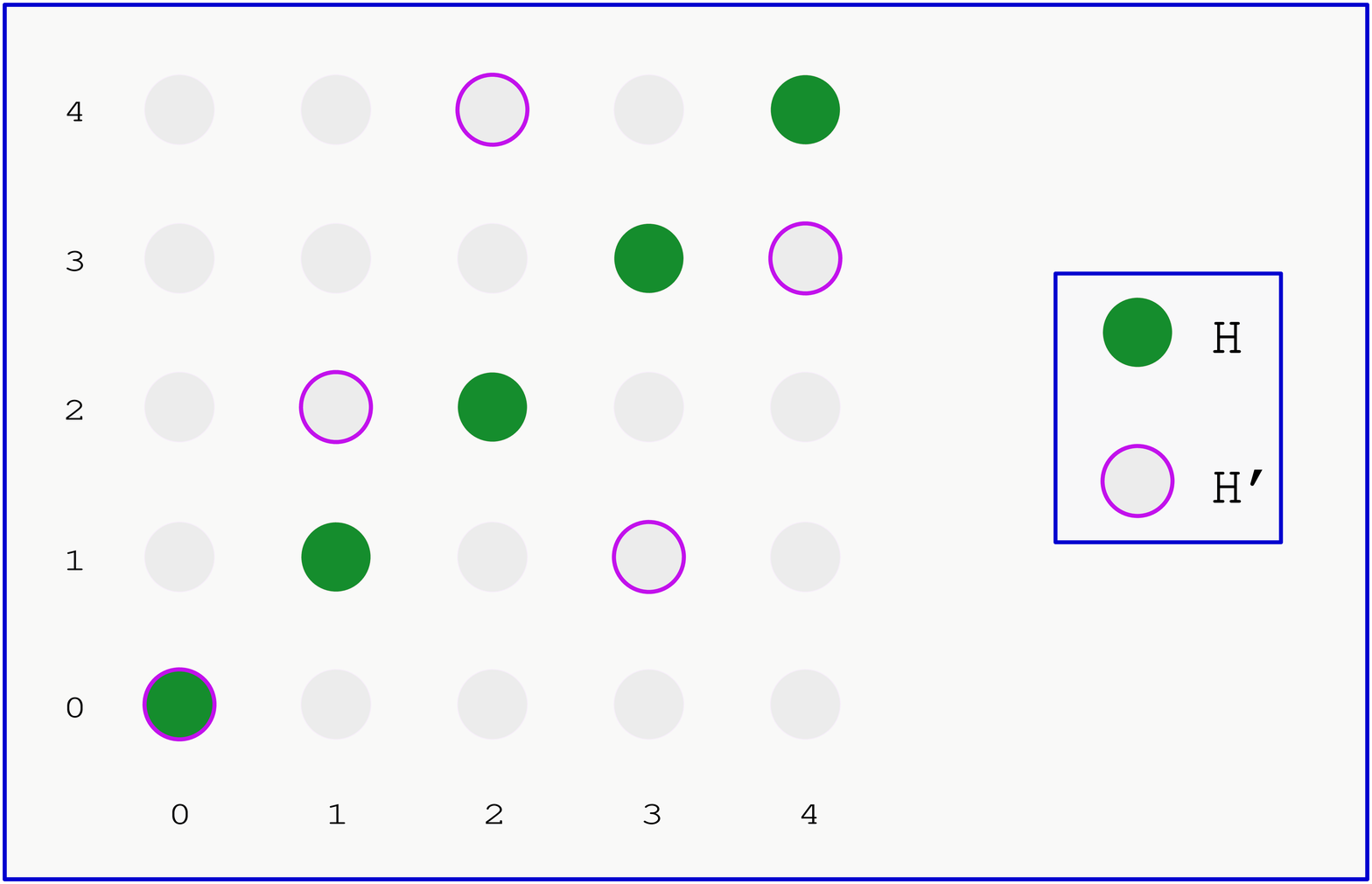}
\end{center}

  {\textsc{Example 3:}}\\ 

Let $l_1$, $l_2$, $r_1$ and $r_2\in F$ so that $\phi(l_1) =
(10,0)$, $\phi(l_2) = (35,0)$, $\phi(r_1) = (0,15)$, and $\phi(r_2) =
(0,20)$.  Further, let $g_1$, $g_2\in F$ so that $\phi(g_1) = (2,6)$ and
$\phi(g_2) = (14,4)$.  Let $H = \langle F', l_1,r_1,g_1\rangle$ and let
$H' = \langle F',l_2,r_2,g_2\rangle$.

It is immediate that $\inner{H} = \fk{10}{15}$, $\outer{H} =
\fk{2}{3}$, $\inner{H'}=\fk{35}{20}$, and $\outer{H'} = \fk{7}{4}$.
Both $H$ and $H'$ have residue $5$.  If we apply $\tau_{(2,3)}$ to
$\outer{H}$ and $\tau_{(7,4)}$ to $\outer{H'}$, and draw our
fundamental $5\times5$ rectangle in $Z^2$, we get the following
diagram.  (Below, we are considering $H$ and $H'$ after the rescaling.)

\begin{center}
\includegraphics[height = 140 pt, width = 200 pt]{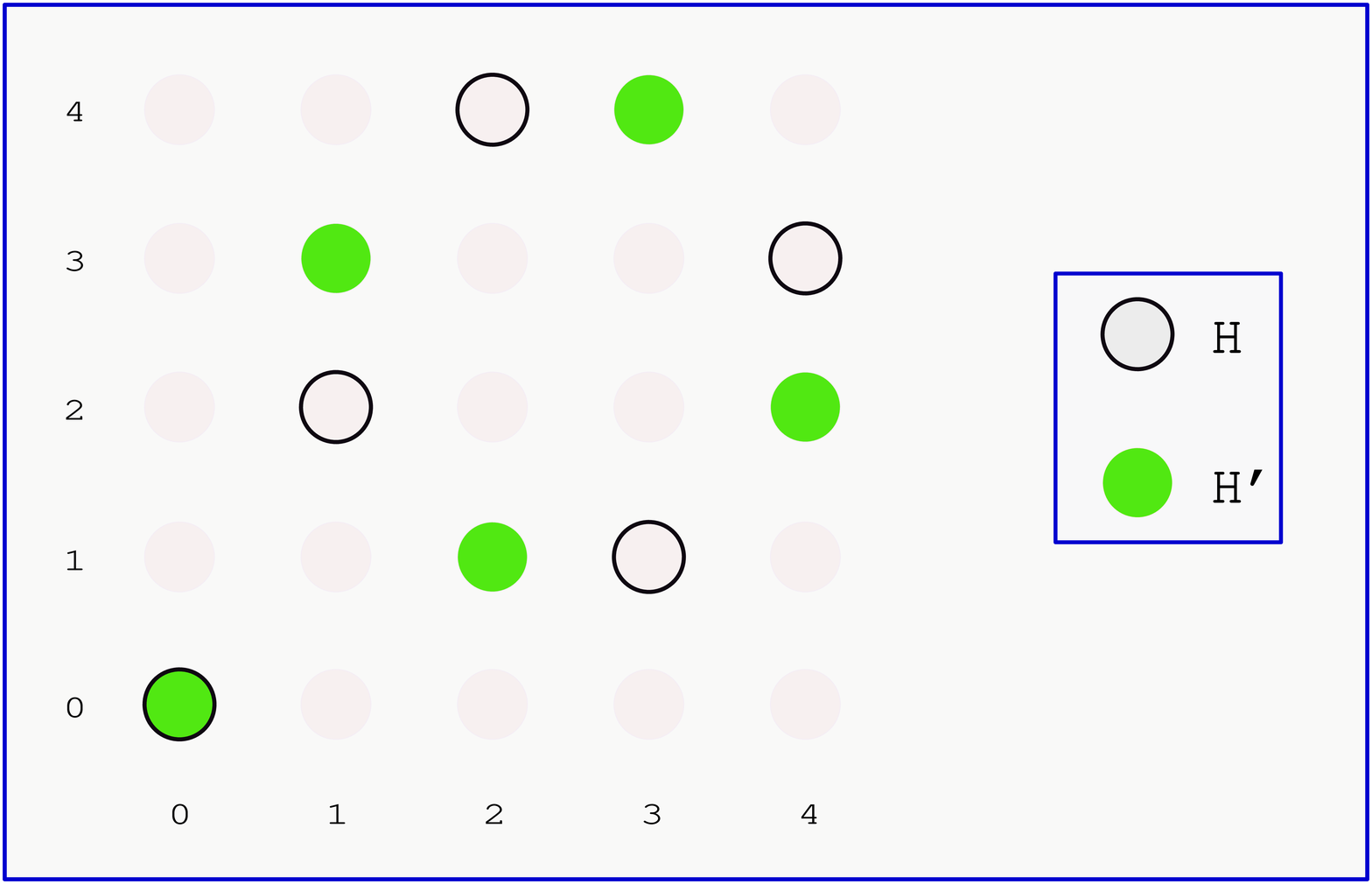}
\end{center}

The scaled version of $H$ is $Rev$ of the scaled version of $H'$, so
that $H\cong H'$.

\newpage
\bibliographystyle{amsplain} \bibliography{ploiBib}

\end{document}